# Complexity Analysis in Bouncing Ball Dynamical System

BY

L. M. Saha[1], Til Prasad Sarma[2] and Purnima Dixit[3]

[1]IIIMIT, Shiv Nadar University, Gautam budh Nagar,U.P.201314, E-mail: lmsaha.msf@gmail.com

[2]Department of Education in Science and Mathematics, NCERT, Sri Aurobindo Marg,
New Delhi -110 016, E-mail: tpsarma@yahoo.com

[3]Department of Mathematics, Dayal Singh College, New Delhi-1100,
E-mail: purnimadixit9@gmail.com

Abstract:

Evolutionary motions in a bouncing ball system consisting of a ball having a free fall in the Earth's gravitational field have been studied systematically. Because of nonlinear form of the equations of motion, evolutions show chaos for certain set of parameters for certain initial conditions. Bifurcation diagram has been drawn to study regular and chaotic behavior. Numerical calculations have been performed to calculate Lyapunov exponents, topological entropies and correlation dimension as measures of complexity. Numerical results are shown through interesting graphics.

1. Introduction

A simple system evolves in simple ways but a complex or complicated system evolve in complicated ways and between simplicity and complexity there cannot be a common ground [1]. Chaos and irregular phenomena may not require very complicated equations. Complexity in a dynamical system can be viewed as its systematic nonlinear properties. It is the order that results from the interaction among multiple agents within the system. A system is complex means its evolutionary behavior do not show regularity but chaotic or some other kind of irregularity. Complexity and chaos observed in a system can well be understood by measuring elements like Lyapunov exponents (LCEs), topological entropies, correlation dimension etc. Topological entropy, a non-negative number, provides a perfect way to measure complexity of a dynamical system. For a system, more topological entropy means the system is more complex. Actually, it measures the exponential growth rate of the number of distinguishable orbits as time advances [2, 3]. Though, positivity measure of Lyapunov exponents (LCEs) signifies presence of chaos, LCEs, topological entropies and correlation dimensions all these three together provide measure of complexities in the system. Motion of a bouncing ball system represented with equations in coupled form of variables, have been appeared in various literatures [4–8] and regular and chaotic motions observed during evolution have been discussed. The models discussed vary with different kind of assumptions and so the variation of nonlinearities.

The present article consisting of a model of bouncing ball system occurring due to a free fall of a ball in the Earth's gravitational field and impacting kinematically certain forced plate [9]. Bifurcation diagrams have been drawn to study some characteristic evolutionary phenomena, (e.g., chaos adding), with increasing numerical value of the driving frequency. Prior to this chaos adding, one observes period doubling bifurcation followed by chaos. Periodic windows appearing within chaos are also subject to study. Numerical investigations carried forward to obtain Lyapunov exponents (LCEs), topological entropies and correlation dimensions for different sets of parameters of the system. Results obtained are shown through graphics.

1. **Bouncing Ball Model**:

Neglecting the air drag, free fall motion of the bouncing ball, with a restitution coefficient k < 1 be written as, [9],

$$\phi_{n+1} = \phi_n + q\, v_n,$$
$$v_{n+1} = k\, v_n + (1+k)\cos(\phi_{n+1}) \qquad (1.1)$$

The system contains an another parameter, q, *the driving frequency* and $q \gg 1$. Bifurcation diagrams for above system are drawn with k = 0.3 and different ranges of values of q shown, respectively, by Fig.1, (Figures (a) - (f)). These shows, initially, the system evolving with period doubling bifurcation followed by chaos and then, a chaos adding phenomena with increasing values of q.

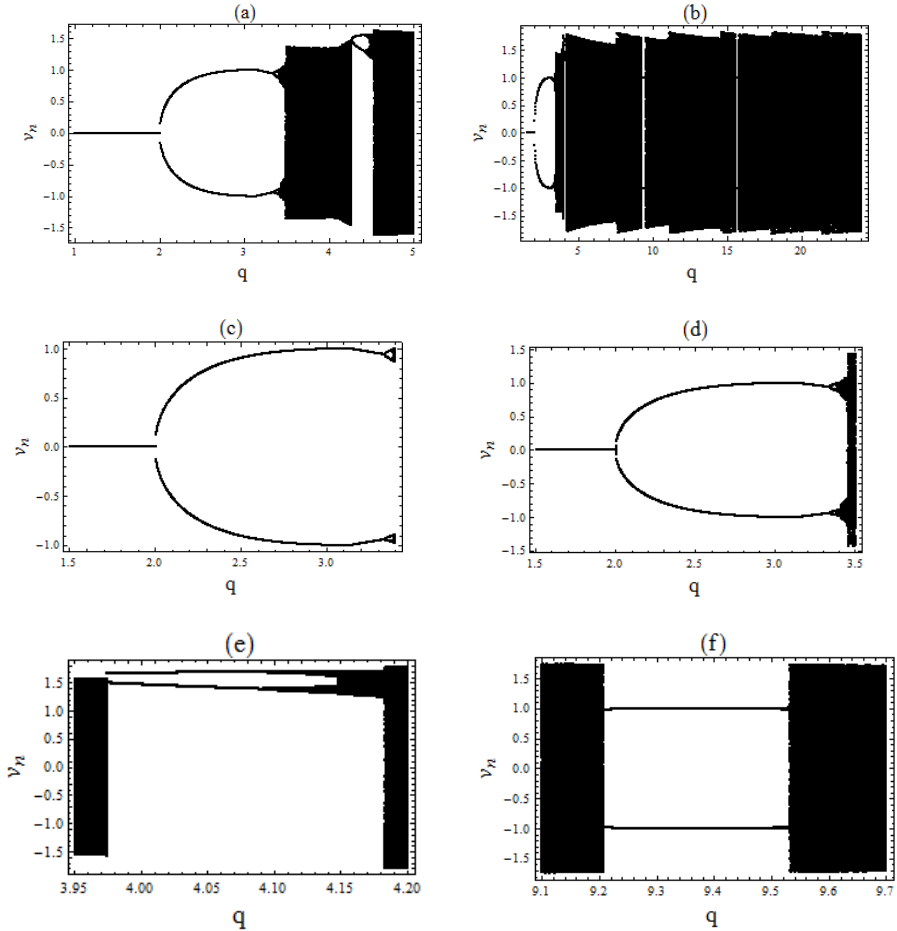

Fig. 1: Bifurcation diagrams of map (1) with k = 0.3 and varying q: (a) $1 \leq q \leq 5$, (b) $1.5 \leq q \leq 24$, (c) $1.5 \leq q \leq 3.4$, (d) $1.5 \leq q \leq 3.5$, (e) $3.95 \leq q \leq 4.2$, (d) $9.1 \leq q \leq 9.7$.

As it appears through bifurcations, before evolving into chaos, system shows regularity for certain rage of parameter value of q while keeping k fixed, k = 0.3. In Fig.2, the figures in upper row are two time series and a two periodic regular attractor for q = 2.5. The lower row figures correspond to those of the upper row for chaotic case when q = 3.8.

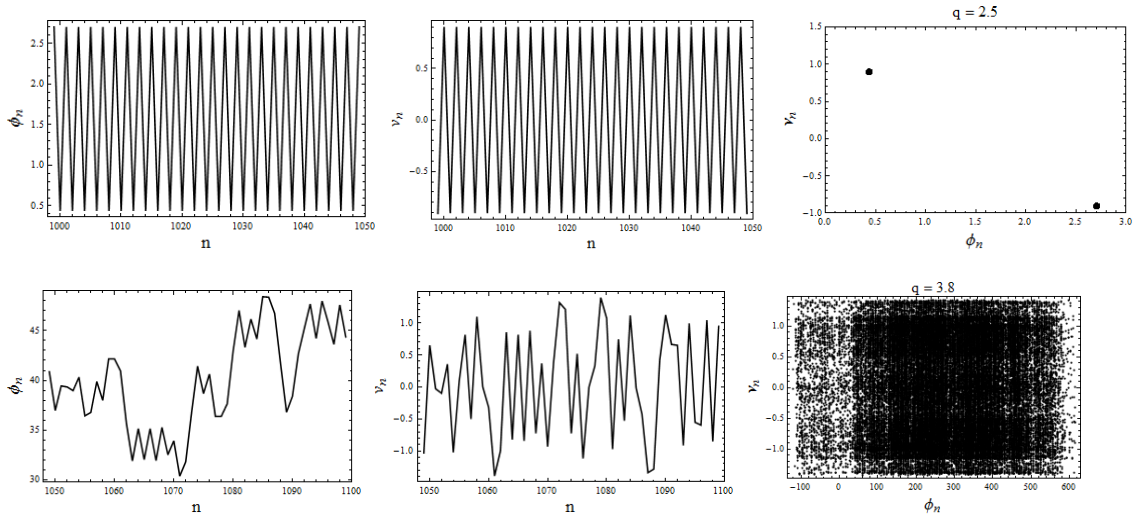

Fig.2: Upper row plots showing clearly 2-periodic attractors for k = 0.3 and q = 2.5; the corresponding plots in the lower row, when q is changed to the value q = 3.8, shows motion is chaotic.

2. **Lyapunov Exponents (LCEs), Topological Entropies & Correlation Dimensions**:

**Lyapunov exponents**: LCEs, have been calculated and plotted, shown in Fig. 3, for regular and chaotic cases as discussed above. In regular case, though the LCEs are negative at each iteration, their numerical values are different. Similarly, for chaotic case, values are pos itive but different.

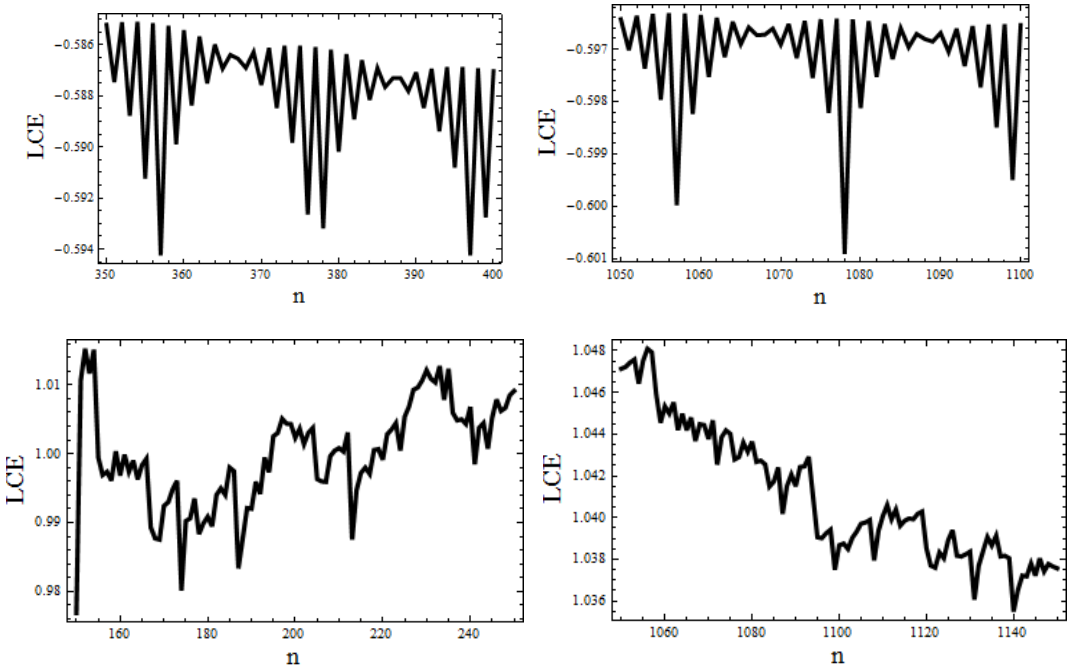

Fig. 3: Plots of LCEs for chaotic and regular cases for fixed k = 0.3. Figures in the upper row are for regular case when q = 3.8 and those of lower row are for q = 2.5.

Magnitude of the positive values of LCEs provide the answer that how chaotic the system be. These differences explain the complexity within the system. With k = 0.3 and q, approximately, q = 3.31, the system shows regular behavior as it is evident from the bifurcation diagrams as well as from the LCEs plot shown in Fig. 4.

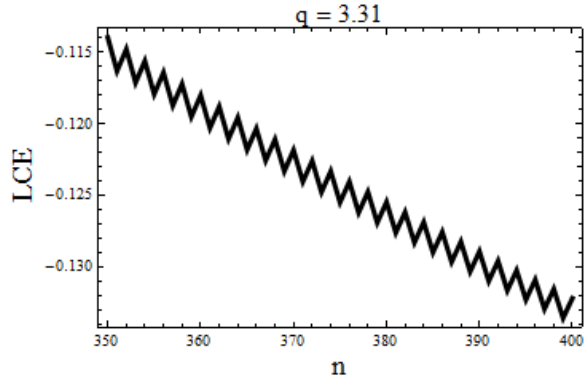

Fig. 4: A plot of LCEs for k = 0.3 and q = 3.31. As LCEs are negative, motion is regular

**Topological Entropy**: Next, let us have calculated topological entropy for the bouncing ball system (1.1) and plotted in Fig. 5.

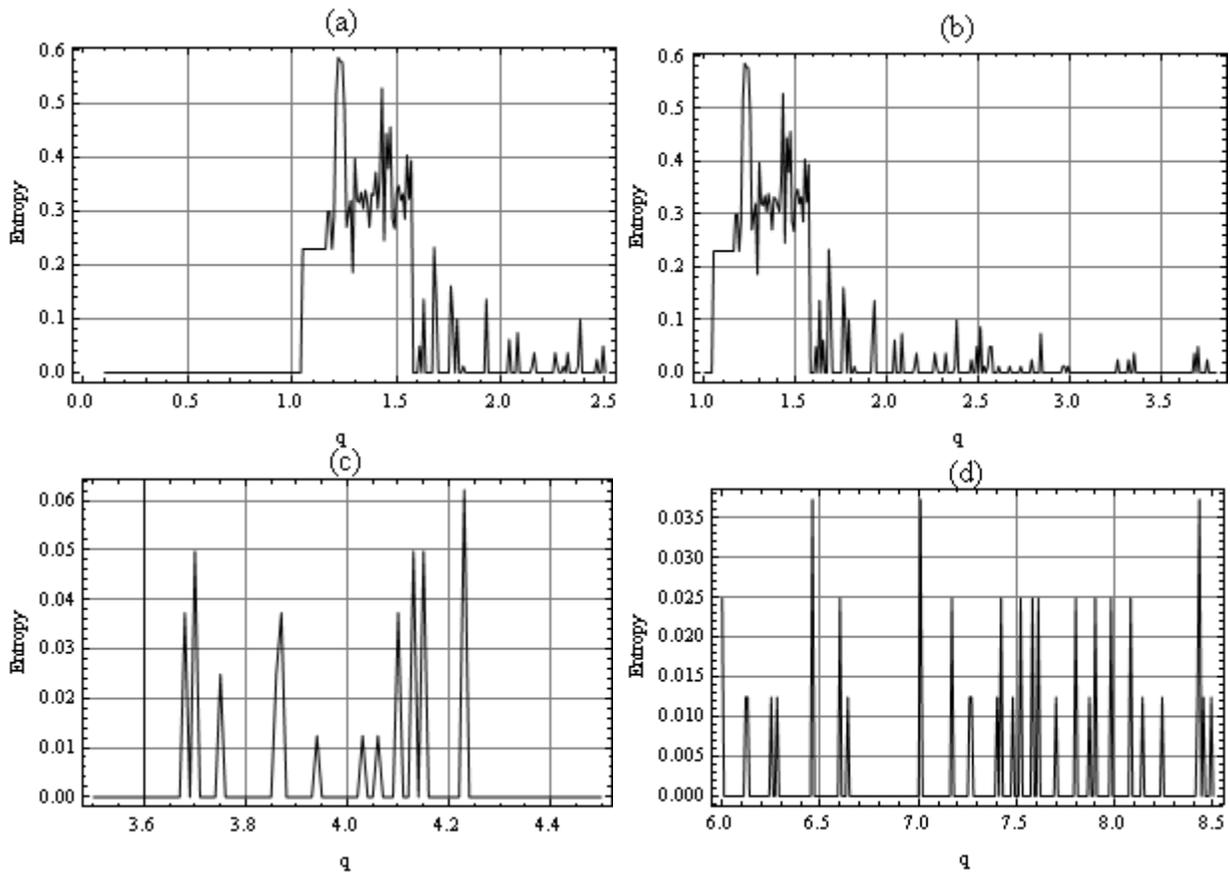

Fig. 5: Plots of topological entropy for k = 0.3 and four ranges of values of parameter q; (a) $0.1 \leq q \leq 2.5$, (b) $1 \leq q \leq 3.8$, (c) $3.5 \leq q \leq 4.5$, (d) $6 \leq q \leq 8.5$.

With k = 0.3 and q approximately, q = 3.31, what we observed in Fig. 4 in LCEs plot, the results obtained here are very different. For $1 \leq q \leq 2.6$, in the former case the system is non-chaotic and shows regularity but in this later case one obtains significant value of topological entropy. As topological entropy measures the complexity, though the system is regular, it is complex. Thus, a non-chaotic nonlinear system can also be complex one. Next, we have plots of 3-Dimension image of topological entropy for $1 \leq q \leq 3.8$, $0.1 \leq k \leq 0.6$, and shown in Fig. 6, which clear picture of complexity.

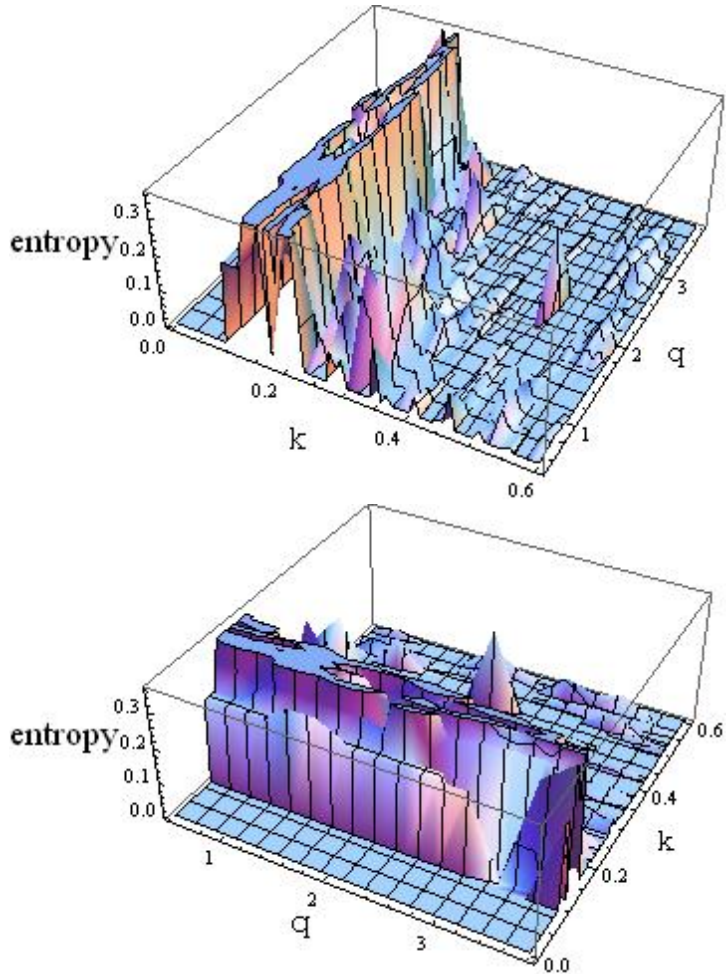

Fig. 6: Two 3-D plots of topological entropies for $1 \leq q \leq 3.8$, $0.1 \leq k \leq 0.6$.

**Correlation dimensions**:

A chaotic set, has fractal structure and so, its correlation dimension gives its measure of dimensionality. Being one of the characteristic invariants of nonlinear system dynamics, the correlation dimension actually gives a measure of complexity for the underlying attractor of the system. To determine correlation dimension a statistical method can be used. It is an efficient and practical method then other methods, like box counting etc. The procedure to obtain correlation dimension follows from some perfect steps calculation [10 – 12]:

Consider an orbit $O(\mathbf{x}_1) = \{\mathbf{x}_1, \mathbf{x}_2, \mathbf{x}_3, \mathbf{x}_4, \ldots\}$, of a map $f: U \to U$, where U is an open bounded set in $\mathbb{R}^n$. To compute correlation dimension of $O(x_1)$, for a given positive real number r, we form the correlation integral,

$$C(r) = \lim_{n \to \infty} \frac{1}{n(n-1)} \sum_{i \neq j}^{n} H\left(r - \|x_i - x_j\|\right), \quad (3.1)$$

where

$$H(x) = \begin{cases} 0, & x < 0 \\ 1, & x \geq 0 \end{cases},$$

is the unit-step function, (Heaviside function). The summation indicates the the number of pairs of vectors closer to r when $1 \leq i, j \leq n$ and $i \neq j$. C(r) measures the density of pair of distinct vectors $\mathbf{x_i}$ and $\mathbf{x_j}$ that are closer to r.

The correlation dimension $D_c$ of $O(\mathbf{x}_1)$ is defined as

$$D_c = \lim_{r \to 0} \frac{\log C(r)}{\log r} \quad (3.2)$$

To obtain $D_c$, log C(r) is plotted against log r and then we find a straight line fitted to this curve. The intercept of this straight line on y-axis provides the value of the correlation dimension $D_c$.

In case of bouncing ball, with k = 0.3 and q = 3.8, the obtained correlation curve is shown in Fig. 7.

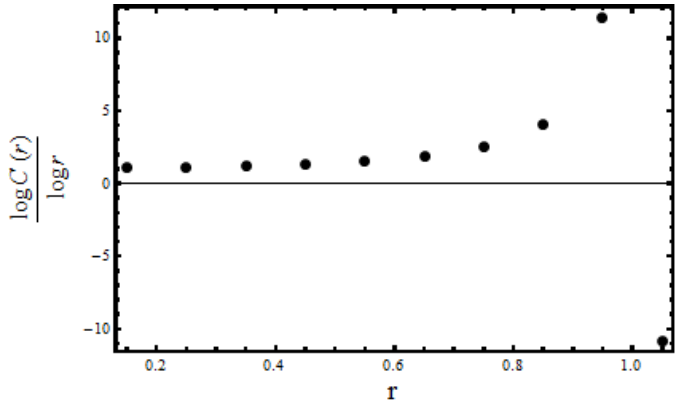

Fig. 7: Plot of correlation curve for bouncing ball system (2.1)

A linear fit to the correlation data be obtained as

$$Y = 2.1792 - 1.0374\ x \quad (3.3)$$

The y – intercept of this straight line is 2.1792 and so the correlation dimension of the chaotic set be measured, approximately, as $D_c = 2.197$.

3. **Discussions**:

The problem of complexity and chaos observed during evolution of bouncing ball has been studied in detail and certain measures of complexity such as Lyapunov exponents, topological entropies, correlation dimension are calculated. The non-negative real number, the topological entropy, describes a perfect measure complexity of dynamical system in the sense that more the topological entropy a system has means it is more complex. Actually, a topological entropy measures the exponential growth rate of the number of distinguishable orbits as time advances in the system. However positivity of its value does not justify the system be chaotic. For k = 0.3 the system studied in this article, chaos happens when q values be increased from 3.31 (approximately). Before this the system is regular. However, we find even much before q reaching this value, within the range $1 \leq q \leq 2.6$, the system is complex as the topological entropies are more than zero. Also, as shown in Fig. 5, for k = 0.3 and values of q in intervals $3.5 \leq q \leq 4.5$ and $6 \leq q \leq 8.5$, where the system showing highly chaotic, topological entropy appears to be very low. Another interesting thing be observed the correlation dimension for q = 2.4, (a regular case), is non-zero and it is given by $D_c$ = 1.653 (approximately). However, there are nonlinear systems

Finally, one can conclude in case of bouncing ball dynamics, complexity and chaos are certainly mixed phenomena.

References:
[1] Ian Stewart. Does God Play Dice ?, Penguin Books 1989

[2] Adler, R L Konheim, A G McAndrew, M H. Topological entropy, Trans. Amer. Math. Soc. 1965; 114: 309-319
[3]R. Bowen, R . Topological entropy for noncompact sets, Trans. Amer. Math. Soc. 1973: 184: 125-136
[4] Holmes PJ. The dynamics of repeated impacts with a sinusoidally vibrating table. J Sound Vib 1982;84:173–89.
[5] Everson RM. Chaotic dynamics of a bouncing ball. Physica D 1986;19:355–83.
[6] Mello TM, Tuffilaro NM. Strange attractors of a bouncing ball. Am J Phys 1987;55:316–20.
[7] Kini A, Vincent TL, Paden B. The bouncing ball apparatus as an experimental tool. J Dyn Syst Meas Control 2006;128:330–40.
[8] Sebastian Vogel, Stefan J. Linz, Regular and chaotic dynamics in bouncing ball models. International Journal of Bifurcation and Chaos, Vol. 21, No. 3, 2011: 869–884
[9] Litak G, Syta A, Budhraja M, Saha, L M. Detection of the chaotic behaviour of a bouncing ball by the 0–1 test. Chaos, Solitons and Fractals 42, 2009: 1511–1517
[10] Grassberger P, Procaccia I. Measuring the Strangeness of Strange Attractors, Physica 9D, 1983:189-208.
[11] Martelli M. Introduction to Discrete Dynamical Systems and Chaos. John Wiley & Sons, Inc., 1999, New York.
[12] Nagashima H, Baba Y . Introduction to Chaos: Physics and Mathematics of Chaotic Phenomena. Overseas Press India Private Limited, 2005.